\documentclass[12pt, reqno]{amsart}

\usepackage[ruled]{algorithm2e}
\usepackage{graphicx}
\usepackage{amsmath}

\newtheorem{thm}{Theorem}
\newtheorem{lemma}[thm]{Lemma}
\newtheorem{prop}[thm]{Proposition}
\newtheorem{cor}[thm]{Corollary}
\theoremstyle{definition}
\newtheorem{defn}[thm]{Definition}

\newtheorem{remark}[thm]{Remark}

\newcommand{\ignore}[1]{}

\newcommand{\CC}{\mathcal{C}}
\newcommand{\MM}{\mathcal{M}}

\begin{document}

\title{The Neighbor-Net Algorithm}
\author{Dan Levy \and Lior Pachter}
\address{Department of Mathematics and Computer Science, UC Berkeley and
the Department of Statistics, University of Oxford}
\email{lpachter@math.berkeley.edu}
\dedicatory{\today}

\begin{abstract}

The neighbor-joining algorithm is a popular phylogenetics method for
constructing trees from dissimilarity maps.
The neighbor-net algorithm is an extension of the neighbor-joining
algorithm and is used for constructing split networks.
We begin by describing the output of neighbor-net in terms of the
tessellation of $\overline{\MM}_{0}^n(\mathbb{R})$ by associahedra. This
highlights the fact that neighbor-net outputs a tree in addition to a
circular ordering and we explain when the neighbor-net tree is the
neighbor-joining tree. A key observation is that the tree constructed in
existing implementations of neighbor-net is not a neighbor-joining tree.
Next, we show that neighbor-net is a greedy algorithm for finding
circular split systems of minimal balanced length. This
leads to an interpretation of neighbor-net as a greedy algorithm for
the traveling salesman problem. The algorithm is optimal for Kalmanson
matrices, from which it follows that neighbor-net is consistent and has
optimal radius $\frac{1}{2}$. We also provide a statistical
interpretation for the balanced length for a
circular split system as the length based on weighted least squares
estimates of the splits.
We conclude with applications of these results and demonstrate the
implications of our theorems for a recently published comparison of
Papuan and Austronesian languages.

\end{abstract}

%\thanks{Supported in part by ....}
%\subjclass{Primary 55M20; Secondary 47H10}
\maketitle
\section{Introduction}

The neighbor-net algorithm was introduced by Bryant and Moulton in
\cite{Bryant2004}. It is a method for  constructing split networks
\cite{Dress2004} from
distance measurements, and has been used for evolutionary analyses
 in  linguistics \cite{Bryant2005b,Dunn2005} and phylogenetics \cite{Huson2005b}.
%The article describing neighbor-net has been cited over 57 times\footnote{source: ISI citation index}.
Neighbor-net is gaining in popularity because it is as fast as
distance based methods for tree construction, and the split networks
output by the algorithm are
informative for studying
conflicting signals in data. The interpretations of
split networks are based on
$T$-theory \cite{Bandelt1992,Dress1996}, which is an active research
area within mathematics.

Despite the intuitive appeal of split networks for data analysis,
a criticism of their use in phylogenetics, and of the neighbor-net algorithm in particular, has been the lack of an obvious tree interpretation. Moreover, although it was remarked in \cite{Bryant2004}
 that ``neighbor-net is based on the
neighbor-joining algorithm of Saitou and Nei
\cite{Saitou1987}'', this was meant to indicate analogy at a high
level: neighbor-net and neighbor-joining are both agglomerative
algorithms, they have similar  selection criteria, and they are
both consistent. However despite the obvious similarities between neighbor-net
and neighbor-joining, there has been no direct link established between the
outputs of the algorithms. It is desirable to establish a mathematically
precise connection because
there have been a number of recent papers ``explaining''
neighbor-joining \cite{Gascuel2006}, both in terms of showing what it
optimizes \cite{Desper2005} and why it works well in practice \cite{Mihaescu2006}. The lack of informative theorems about neighbor-net coupled with the
difficulties in mastering $T$-theory have contributed to
a sense that  interpretations of neighbor-net results ``remain messy and
subject to a certain degree of subjectivity\footnote{The statement
appears in the specific context of a commentary on a paper describing the
classification of Bantu languages \cite{Lutz2006}; we believe
that it reflects prevailing sentiment about the neighbor-net algorithm
and its utility for evolutionary analyses.}.''

We describe the precise connection between neighbor-net and
neighbor-joining in Section 2, and in Section 5 we show that our observation can be used to
allay concerns  that neighbor-net provides no direct
phylogenetic {\em tree} information. Our result also provides
an interpretation of $\overline{\MM}_{0}^n(\mathbb{R})$ as the
  space of phylogenetic networks.
In Section 3 we show that neighbor-net is a
greedy algorithm for the traveling salesman problem that minimizes the
balanced length of the split system at every step. This extends the
notion of balanced length in \cite{Semple2004} and the results of
\cite{Desper2005} where it was shown that neighbor-joining greedily optimizes the
balanced length of a tree. In Section 4, we prove that neighbor-net is optimal for Kalmanson
dissimilarity maps. This establishes new proofs for results of
\cite{Christopher1997,Christopher1996,Deineko1998}, and provides an
analog of Atteson's neighbor-joining robustness theorem \cite{Atteson1999} for neighbor-net.

\section{The mathematics}

The main objects of study in this paper are a class of discrete
metric spaces called {\em circular decomposable metrics} that include
tree metrics as a special case. We begin with an introduction to some
fundamental results about these metric spaces. Their
study is part of $T$-theory, and we refer the reader to \cite{Dress1996} for a
more thorough introduction and survey of the subject.
Throughout the
paper, $X=\{1,\ldots,n\}$ denotes the finite set on which metrics are defined.
\begin{defn}
A {\em split} $S=\{A,B\}$ is a partition of $X$ into two non-empty blocks.
A set of splits is called a {\em split system}.
The split metric determined by $S$ is the pseudo-metric
\[ \delta_S =
\begin{cases}
0 & \mbox{if } \{x,y\} \subseteq A \mbox{ or } \{x,y\} \subseteq B.\\
1 & \mbox{otherwise}.
\end{cases}\]
\end{defn}
\begin{defn}
A split system $\mathcal{S}$ is {\em pairwise compatible} if
 for every pair of distinct splits
  $S_1=\{A,B\}$, $S_2=\{A',B'\}$ in $\mathcal{S}$, at least one of the
  intersections \[ A \cap A',\, A \cap B',\, B \cap A',\, B \cap B'\] is empty.
\end{defn}
\begin{defn}
A {\em dissimilarity map} on $X=\{1,\ldots,n\}$ is a function $\delta:X \times X \rightarrow
  {\mathbb R}$ that satisfies $\delta(i,j)=\delta(j,i) \geq 0$
  and $\delta(i,i)=0$.
A dissimilarity map $\delta$ satisfies the {\em four point condition} if
for every four elements $i,j,k,l \in X$, two of the three terms in the
following list are equal and greater than the third:
\[ \delta(i,j) + \delta(k,l), \, \delta(i,k)+\delta(j,l), \, \delta(i,l)+\delta(j,k).\]
\end{defn}
\begin{thm}[\cite{Semple2003}]
\label{thm:fourpoint}
The following are equivalent statements about $\delta:X \times X
\rightarrow \mathbb{R}$:
\begin{enumerate}
\item There exists a split system $\mathcal{S}$ such that every pair of
  distinct splits in $\mathcal{S}$ is pairwise compatible, and $\delta = \sum_{S \in \mathcal{S}} \lambda_S \delta_S$ where $\lambda_S \geq 0$ for
  all $S \in \mathcal{S}$.
\item $\delta$ is a metric and satisfies the four point condition.
\end{enumerate}
\end{thm}
There is a canonical median graph associated with a split system called
the {\em Buneman graph} \cite{Buneman1971}. The Buneman graph of a
pairwise compatible split system is a tree, and therefore, in light
of Theorem \ref{thm:fourpoint}, metrics satisfying the four point condition are
called {\em tree metrics}. They are precisely the metrics $\delta:X
\times X \rightarrow \mathbb{R}$ for
which there is an edge weighted tree whose leaves are labeled by $X$, and for which
$\delta(i,j)$ is the ``additive distance'' between $i$ and $j$ in the tree.

Theorem \ref{thm:fourpoint} provides the necessary ingredients for
describing the input and output of the neighbor-joining
algorithm. Specifically, neighbor-joining is an efficient algorithm for
evaluating a certain function from the set of
dissimilarity maps to pairwise compatible split systems. A key feature of the algorithm, is
that the steps explicitly construct the Buneman tree associated with the
output.

The neighbor-net algorithm is similarly explained in terms of certain
split systems and metrics. The key concept is that of a circular
ordering for a finite set $X$.
\begin{defn}
A {\em circular ordering} $\pi = \{x_1,\ldots,x_n\}$ is a bijection
between $X$ and the vertices of the $n$-cycle $C_n$ such
that $x_i$ and $x_{i+1}$ are adjacent vertices of $C_n$. We adopt the
convention that $x_{n+1}=x_1$.
\end{defn}
Given a circular ordering $\pi$, let $W_{\pi} = \{ \{\{x_i,x_j\},\{x_k,x_l\}\}:i<j<k<l \mbox{ or } l<i<j<k\}$. Note that $W_{\pi}$ is a set consisting of pairs of sets constructed from quartets. In what follows we use the notation $(ij;kl)$ to denote the quartet $\{\{x_i,x_j\},\{x_k,x_l\}\}$.
\begin{defn}
A split system $\mathcal{S}$ is {\em circular} with respect to a circular
ordering $\pi=\{x_1,\ldots,x_n\}$ if every split $S \in \mathcal{S}$ is
of the form
\[ S = \{ \{x_{i+1},\ldots,x_j\},\{x_{j+1},\ldots,x_i \}  \} \mbox{ for some } i<j. \]
\end{defn}
Note that every pairwise compatible split system is circular.
\begin{defn}
A dissimilarity map $\delta$ satisfies the
{\em Kalmanson conditions} \cite{Kalmanson1975} with respect to a circular ordering
$\pi$ if for every $i < j < k < l$,
\begin{eqnarray*}
\delta(x_i,x_j) + \delta(x_k,x_l) & \leq & \delta(x_i,x_k) + \delta(x_j,x_l), \\
\delta(x_i,x_l) + \delta(x_j,x_k) & \leq & \delta(x_i,x_k) + \delta(x_j,x_l).
\end{eqnarray*}
\end{defn}
Given a dissimilarity map $\delta$ that satisfies the Kalmanson conditions with respect to a circular ordering $\pi$, we let $W_{\delta} = \{(ij;kl):\delta(x_i,x_j)+\delta(x_k,x_l)< \delta(x_i,x_k)+\delta(x_j,x_l) \mbox{ for } i<j<k<l \mbox{ or } l<i<j<k  \}$. Note that $W_{\delta} \subseteq W_{\pi}$ is a set of quartets given by the strict Kalmanson inequalities.

%We say that a quartet $(x_i,x_j; x_k, x_l)$ satisfies a {\em strict Kalmanson condition} iff
%\begin{eqnarray*}
%\delta(x_i,x_j) + \delta(x_k,x_l) & < & \delta(x_i,x_k) + \delta(x_j,x_l).
%\end{eqnarray*}

\begin{thm}[\cite{Chepoi1998,Christopher1996}]
\label{thm:circular}
The following are equivalent statements about $\delta:X \times X
\rightarrow \mathbb{R}$:
\begin{enumerate}
\item There exists a circular ordering $\pi$ and a split system
  $\mathcal{S}$ so that
$\delta = \sum_{S \in \mathcal{S}}\lambda_S\delta_S$ where
every split $S \in \mathcal{S}$ is circular with respect to $\pi$ and
  $\lambda_S \geq 0$ for all $S \in \mathcal{S}$.
\item $\delta$ is a metric and satisfies the Kalmanson conditions with respect to $\pi$.
\end{enumerate}
Moreover, a quartet $(ij;kl) \in W_{\delta}$ iff there exists a split $S$ with $\lambda_S>0$ such that $i,j$ and $k,l$ are in different blocks of $S$.

%Further, $\delta$ satisfies the strict Kalmanson condition for all quartets $(a_1,a_2; b_1, b_2)$
%where $a_1, a_2 \in A$ and $b_1, b_2 \in B$ if and only if $\lambda_S > 0$ for $S = {A, B}$.
\end{thm}

Metrics satisfying condition (1) of Theorem \ref{thm:circular} are
called {\em circular decomposable metrics}, and it is possible to represent them
using {\em split graphs}. These are described in
detail in \cite{Bryant2004}. Here we merely illustrate the idea
with an example (Figure 1(a,b)). Each class of parallel edges
corresponds to one split $S \in \mathcal{S}$ and the length of the edges
in a class are given by the $\lambda_S$. Split graphs are not necessarily
unique, but they provide a useful way to visualize a circular
decomposable metric.
%\begin{remark}
%If $S \in \mathcal{S}$ is a split of a circular split system with $\lambda_S >0$, then for any $i,j$ and $k,l$ in different blocks of $S$ we have that $\{\{i,j\},\{k,l\}\} \in Q_{\delta}$.
%\end{remark}
The neighbor-net algorithm  outputs a circular ordering for the
purpose of visualizing a circular decomposable metric associated to
it using split graphs. The algorithm is {\em agglomerative}, which means that the circular
ordering is constructed iteratively. The boxed Algorithm 1 describes the
details of the algorithm. The terms used in its description are defined
below:
\begin{algorithm}
\label{neighbornet}
  \SetLine
  \AlgData{A dissimilarity map $\delta:X \times X \rightarrow \mathbb{R}$.}
  \AlgResult{Circular ordering $\pi:X \rightarrow C_n$ together with a
 split system $\mathcal{T}$ of $n-1$ pairwise compatible splits that are
    circular with respect to $\pi$.}
  Let $G$ be the disjoint union of $n$ vertices and $\CC$ the partial
  circular ordering with graph $G$. Let $\mu:X \rightarrow
  \mathbb{R}$ be the weighting for $\CC$. \\
  \While{$|\CC|>1$}{
    \For{$i,j \in {|\CC| \choose 2}$}{
      Set $Q_{\delta}(C_r,C_s)=(|\CC|-2)\delta(C_r,C_s)-\sum_{t \in \CC \setminus \{C_r\}}\delta(C_r,C_t) - \sum_{C_t \in \CC \setminus \{C_s\}}\delta(C_t,C_s)$.
    }
[{\bf Selection step part 1}] Choose a pair $C_{r^*},C_{s^*} \in \CC$ that minimizes $Q_\delta$\;
    \For{$i \in \hat{C}_{r^*}, \, j \in \hat{C}_{s^*}$}{
      Set $\hat{Q}_{\delta}(i,j) = (|\CC|-4+|\hat{C_{r^*}}|+|\hat{C_{s^*}}|)\delta(i,j) - \sum_{t \neq r^*,s^*}
\delta(i,C_t) - \sum_{t \neq r^*,s^*} \delta(j,C_t) - \sum_{k \in (C_{r^*} \cup C_{s^*}) \setminus\{i\}} \delta(i,k) -\sum_{k \in (C_{r^*} \cup C_{s^*}) \setminus\{j\}} \delta(j,k)$.}
[{\bf Selection step part 2}] Choose the pair $i^* \in \hat{C}_{r^*}, \, j^* \in \hat{C}_{s^*}$ that minimizes $\hat{Q}_\delta$\;
[{\bf Merge step}]   Let $u,v$ be the vertices in the circular ordering graph corresponding to $i^*$ and $j^*$.  Add the edge $(u,v)$ to the circular ordering graph and coarsen the
  partition $\CC$ by merging $C_{r^*}$ and $C_{s^*}$.

[{\bf Adjustment step}]  Adjust $\mu(i), i \in C_{r^*} \cup C_{s^*}$ so that $\sum_{i \in C_{r^*}\cup C_{s^*}} \mu(i)=1$.

[{\bf Tree construction step}]   Add the split $\{\{C_{r^*} \cup C_{s^*}\},\{\cup_{t \neq r^*,s^*} C_t \} \}$ to the
  distinguished list.
}
Output the circular ordering $\pi$ and the split system $\mathcal{T}$.
  \caption{Neighbor-net algorithm}
\end{algorithm}

\begin{defn}
Let $G$ be a subgraph of the cycle $C_n$ with $n$ vertices and $m$
components.
The graph $G$ is called the {\em circular ordering graph}.
A {\em partial circular ordering} $\CC$ consists of the graph $G$
together with a bijection between
$X$ and the vertices of $G$.

Equivalently, a partial circular ordering is a partition $\CC$ of $X$
into ordered sets
$\CC=\{C_1,\ldots,C_m\}$ where each $C_r \subseteq X$ and
$i,j$ are adjacent elements in $C_r$ for some $r$ iff $i,j$ correspond to adjacent vertices in
$G$. We use the notation $\hat{C}_r$ to denote the vertices of degree 0
or 1 in the subgraph corresponding to $C_r$.
\end{defn}
\begin{defn}
Let $\CC$ be a partial circular ordering with $|\CC|=m$. A {\em weighting}
for $\CC$ consists of a function $\mu:X \rightarrow \mathbb{R}$ such
that $\mu(i) \geq 0$ for all $i \in X$, and for each $r \in \{1,\ldots,m\}$, $\sum_{i \in C_r}, \mu(i) = 1$ and $\mu(i)>0$ for all $i \in \hat{C}_r$.
We define
\begin{eqnarray}
\label{eq:mu}
 \delta(C_r,C_s) &:=& \sum_{i \in C_r, j \in C_s}
\mu(i)\mu(j)\delta(i,j), \mbox{and }\\
 \delta(x, C_r) &:=& \sum_{i \in C_r} \mu(i)\delta(x,i)
\end{eqnarray}
\end{defn}

Note that if $|\CC|=|X|$ then there is only one weighting for $\CC$,
i.e., $\mu(i)=1$ for all $i$.
Next, we introduce two types of weightings that lead to interesting
neighbor-net algorithms in Sections 3 and 5.
\begin{defn}
A weighting
 $\mu:X
\rightarrow \mathbb{R}$ is a {\em TSP weighting} if, for all $i \in X$, $\mu(i) = 0$ for
all $i \notin \hat{C_r}$.
\end{defn}
These weightings lead to aggressive greedy algorithms for the traveling salesman
problem (Theorem \ref{thm:NN_Greedy}).
\begin{defn}
\label{def:tree_weighting}
Let  $\mu:X
\rightarrow \mathbb{R}$ be a weighting for a partial circular ordering
$\CC$, and consider a new weighting $\mu':X \rightarrow \mathbb{R}$
for the adjustment step of neighbor-net. $\mu'$ is a {\em tree
  weighting} if it satisfies
\[  \mu'(i)  =
\begin{cases}
  \alpha \mu(i) & \mbox{ if } i \in C_r,\\
(1-\alpha)\mu(i) & \mbox{ if } i \in C_s,
\end{cases}\]
where $C_r$ and $C_s$ are the two blocks being merged in the merging
step and $0 \leq \alpha \leq 1$.
\end{defn}
Tree weightings are so named because of the following proposition:

\begin{prop}
\label{prop:tree_construct}
The split system $\mathcal{S}$ output by neighbor-net on input $\delta$ is pairwise
compatible, and in bijection with a binary tree $T$.
If $\mu$ is a tree weighting then the tree $T$ is the neighbor
joining tree for $\delta$, where the agglomeration parameter at every
step is given by the tree weighting parameter $\alpha$.
\end{prop}
{\bf Proof}: Note that the addition of an edge to the graph $G$ during
a run of the algorithm results in a coarsening of the partition $\CC$, where two
blocks are merged into one. For this reason, if $S_1=\{A_1,B_1\}$ is a split added
before $S_2=\{A_2,B_2\}$ to $\mathcal{S}$, then either $A_1 \cap A_2 =
\emptyset$ or $A_1 \cap B_2 = \emptyset$. To see that the tree
determined by $\mathcal{T}$ is the neighbor-joining tree, it suffices to note
that selection step 1, together with the adjustment step specified by a
tree weighting, is identical to the agglomeration procedure of
neighbor-joining.  With a tree weighting, selection step 2 and the fixed ordering within
clusters has no effect on the adjustment or tree construction steps.  If we simply omit
the selection step 2 and the merge step, the neighbor-net algorithm reduces to neighbor-joining.\qed

Proposition \ref{prop:tree_construct} justifies the term {\em tree construction
  step} in the neighbor-net algorithm and shows that the output of
neighbor-net is not only a circular ordering,
but also a tree. The connection to the neighbor-joining tree is explored
further in Section 5.

The coarsenings of the partition $\CC$ in the merge
step are also closely related to {\em graph tubings} \cite{Devadoss1999}:

\begin{defn}
Let $G$ be a finite graph. A {\em tube} is a proper nonempty set of
vertices whose induced graph is a proper, connected subgraph of $G$. A
pair of tubes $r,s$ are {\em nested} if $r \subset s$ or $s
\subset r$. They {\em intersect} if they are not nested and $r \cap
s \neq \emptyset$, and two tubes are {\em adjacent} if $r \cap s =
\emptyset$ and $r \cup s$ is a tube. Two tubes are {\em compatible} if
they do not intersect and are not adjacent. A {\em tubing} of $G$ is a set of
tubes that are pairwise compatible.
\end{defn}

\begin{prop}
\label{prop:tubing}
Let $P_{n-1}$ be the path on $n-1$ vertices. A labeling of $P_{n-1}$ is
a bijection from $\{1,\ldots,n-1\}$ to $P_{n-1}$.
The output of neighbor-net is a labeling of $P_{n-1}$ together with a
maximal tubing of its line graph $L(P_{n-1})$.
\end{prop}
{\bf Proof}: Each coarsening of $\CC$ corresponds to a tube in $L(P_{n-1})$.

\begin{figure}
\label{fig:example1}
\begin{picture}(500,500)
\put(50,0){\includegraphics[scale=0.6]{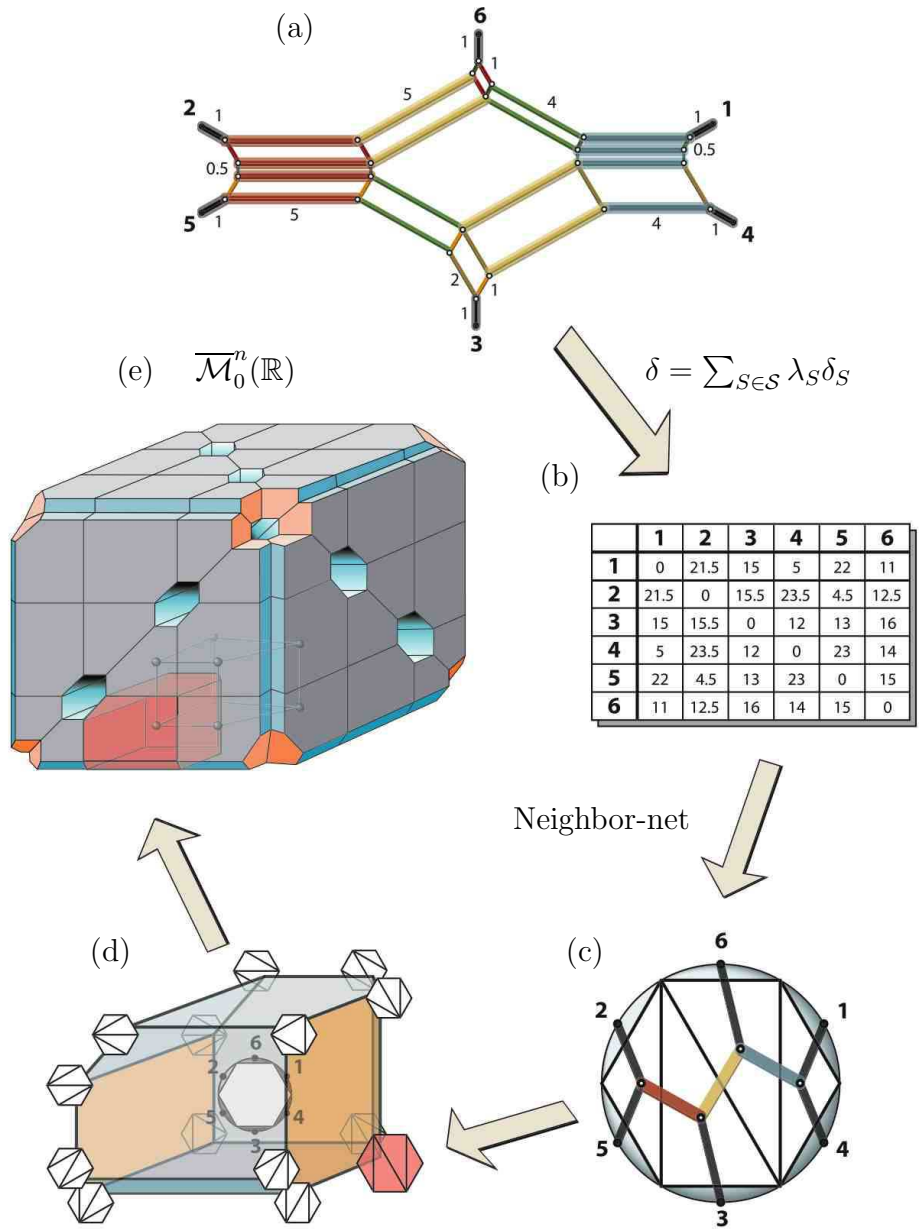}}
\put(160,470){(a)}
\put(300,340){$\delta=\sum_{S \in \mathcal{S}}\lambda_S \delta_S$}
\put(250,170){Neighbor-net}
\put(130,340){$\overline{\MM}_{0}^n(\mathbb{R})$}
\put(260,300){(b)}
\put(100,340){(e)}
\put(90,120){(d)}
\put(270,120){(c)}
\end{picture}
\caption{(a) A split network representation of a circularly decomposable
  metric. Each split $S$ corresponds to a color class with the length of the
  edges in the class $\lambda_S$ indicating the size of the split. (b)
  The metric $\delta$ derived from the splits network. (c) The output
  of  neighbor-net on input $\delta$. The tree is the neighbor-joining
  tree. Note that its edges are  highlighted in the splits
  network. (d) The associahedron $K_5$ corresponding to the circular
  ordering $\pi=\{1,4,3,5,2,6\}$ and the vertex
  corresponding to the neighbor-joining tree. (e) The space of phylogenetic
  networks $\overline{\MM}_{0}^n(\mathbb{R})$.}
\end{figure}

\begin{defn}[\cite{Carr2006}]
For a graph $G$ with $n$ vertices, the graph-associahedron
$\mathcal{P}G$ is the convex polytope of dimension $n-1$ whose face
poset is isomorphic to the set of valid tubings of $G$, with the poset
order corresponding to nesting of tubes.
\end{defn}
The {\em associahedron} (denoted by $K_n$) refers to the graph-associahedron of the
path $P_{n-1}$, and its vertices are in bijection with tubings of the path.
\begin{prop}[See Figure 1(c,d)]
\label{prop:Catalan}
The number of vertices of $K_{n-1}$ is given by
the Catalan number $\frac{1}{n-1}{2n-4 \choose n-2}$. The vertices are
in bijection with tubings of the path $P_{n-2}$, triangulations of the
convex $n$-gon, and rooted binary trees with $n-1$ leaves.
\end{prop}
We have listed just a few of the objects in bijection with the vertices
of $K_n$. In fact, there are dozens of combinatorial objects enumerated
by the Catalan numbers (see \cite{Stanley1999}. In the context of the
neighbor-joining algorithm, Proposition \ref{prop:Catalan} appears as
Proposition 3.1(ii) in \cite{Semple2004}.

Proposition \ref{prop:Catalan} allows us to enumerate the total number
of possible outputs of the neighbor-net algorithm.
\begin{prop}
\label{prop:count}
The number of possible outputs of neighbor-net for n taxa is
\[ \frac{(2n-5)!}{(n-3)!}.\]
\end{prop}
{\bf Proof}: The number of distinct circular orderings (where two
orderings are equivalent under the action of the dihedral group) is
$ \frac{1}{2}(n-1)!$ so the total number of possible outputs is
\begin{equation}
 \frac{1}{n-1}{2n-4 \choose n-2} \cdot \frac{1}{2}(n-1)!
 =  \frac{(2n-5)!}{(n-3)!}.  \end{equation}  \qed

 The first numbers are $1,1,1,6,60,840,15120,332640,8648640,259459200,\ldots$
These numbers also appear in another context in computational biology;
 in genome assembly they
are the number of ways that
$n$ distinguishable equal-length clones can be interleaved to form one
island \cite{Newberg1996}.

Propositions \ref{prop:tubing} and \ref{prop:Catalan} together establish
that the output of neighbor-net is a circular ordering together with
the vertex of an associahedron. Equivalently, it is a labeled convex
$n$-gon together with a triangulation. Thus, it is natural to consider
$\frac{1}{2}(n-1)!$ associahedra corresponding to the distinct circular
orderings.
These associahedra can be glued together in a natural way so
  that faces are identified when the associated subdivisions of the
  $n$-gon differ by twists along the diagonal \cite{Devadoss1999}.
This identification corresponds exactly to the tessellation of
a certain space known as $\overline{\MM}_{0}^n(\mathbb{R})$ by
associahedra.
The space $\overline{\MM}_{0}^n(\mathbb{R})$ consists of the real points
of the
Deligne-Knudsen-Mumford compactification of the moduli space $\MM_{0}^n$
of Riemannian spheres with $n$ labeled punctures. Its tessellation by
associahedra is described in \cite{Devadoss1999}. Figure 1(e) shows the
example for $n=6$. One element from the dual tessellation by
$n-3=3$-dimensional cubes is also shown. Each cube is divided into $8$
octants, and these octants are in bijection with the possible outputs of
neighbor-net (by Proposition \ref{prop:count} there are 840 of them). This is summarized as follows:
\begin{remark}
Neighbor-net is an efficient evaluation of a function from
dissimilarity maps to octants in the dual tessellation by cubes of $\overline{\MM}_{0}^n(\mathbb{R})$.
The vertices of the cube (or equivalently, each associahedron) can be
interpreted as providing the basis
for circular decomposable metrics (networks) together with tubings of the path
that are in bijection with trees (phylogenies). We
therefore refer to $\overline{\MM}_{0}^n(\mathbb{R})$ (or its dual
tiling) as the {\em space of phylogenetic networks}\footnote{The term {\em phylogenetic network}
is also used to denote other objects, e.g. see \cite{Moret2004}.}.
\end{remark}
We note that the relevance
of $\overline{\MM}_{0}^n(\mathbb{R})$ to phylogenetics was already mentioned
in \cite{Billera2001}, however in that paper it was deemed unsuitable for
describing the space of trees, and
 replaced with a quotient space equivalent to the tropical
Grassmanian \cite{ASCB2005}. It is interesting that
$\overline{\MM}_{0}^n(\mathbb{R})$ also appears in the study of genome
rearrangements \cite{Barad2003}.
It should be interesting
to  explore
extensions of neighbor-net that produce, via agglomeration, tubings of
line graphs other than $P_{n-1}$, thus leading to more general
phylogenetic networks connected to graph associahedra.

We conclude this section by noting that our description of neighbor-net
has been based on an interpretation of the algorithm as producing only
combinatorial output, i.e., a circular ordering $\pi$ together with a tree. In
practice, it is possible to obtain weights $\lambda_S$ for the splits in
the circular split system $\mathcal{S}$ compatible with $\pi$ in the course of
the algorithm. This is done by setting
\begin{equation}
\label{eq:distance_estimates}
\lambda_S =
\frac{1}{2}\left(\delta(x_i,x_j)+\delta(x_{i-1},x_{j-1})-\delta(x_i,x_{j-1})-\delta(x_{i-1},x_j)\right)
.\end{equation}
for every split $S=\{\{x_i,\ldots,x_j\},\{x_{j+1},\ldots,x_{i-1}\}\}$.

The problem with such a procedure is that there is no guarantee that all
the $\lambda_S$ will be non-negative, and therefore the result may not
be a circular decomposable metric. This may be circumvented by setting
$\lambda_S$ to zero if it is negative, but this solution may lead to
inaccurate results. For these reasons, a preferable procedure is to use the
circular ordering $\pi$ to subsequently estimate the split weights using
a non-negative least squares optimization method. This was done in the original neighbor-net
implementation \cite{Bryant2004}. 

\section{The computer science}

In the previous section we have explained the input and output of
the neighbor-net algorithm. In this section, we show that
neighbor-net is a greedy algorithm for minimizing the (suitably defined)
length of a dissimilarity map with respect to a circular ordering.
We
begin by extending the formulation of balanced length in
\cite{Semple2004} from trees to circular decomposable metrics.

We say that a circular ordering $\pi=\{x_1,\ldots,x_n\}$ is {\em consistent with
$\CC$}, if for every pair of adjacent elements $i,j$ in some $C_l \in \CC$
there exists a $k$ such that $x_k=i$ and $x_{k+1}=j$.
We denote the
circular orderings consistent with $\CC$ by $o(\CC)$.

\begin{defn}
\label{def:length_pco}
The {\em balanced length of a dissimilarity map} $\delta$ with respect to
a partial circular ordering $\CC$ is defined to be
\begin{eqnarray*}
l(\delta,\CC) & := & \frac{1}{|o(\CC)|}\sum_{(x_1,\ldots,x_n)\in o(\CC)} \left[
    \frac{1}{2} \sum_{i=1}^{n}\delta(x_i,x_{i+1}) \right].\\
 & &  = \frac{1}{2|o(\CC)|}\sum_{(i,j) \in X} \eta_{\CC}(i,j)\delta(i,j).
\end{eqnarray*}
Here $\eta_{\CC}(i,j)$ is the number of circular orderings consistent with $\CC$ where
$i$ is adjacent $j$.
\end{defn}

\begin{remark}
\label{rem:TSP}
The partial circular ordering $\CC^* = {\rm argmin}_{|\CC|=1}(l(\delta,\CC))$ is just
the shortest traveling salesman tour for the dissimilarity map $\delta$.
\end{remark}

We extend the notion of a balanced agglomeration scheme from neighbor
joining to neighbor-net:
\begin{defn}
A {\em balanced TSP weighting} is a TSP weighting where
\[ \mu(i) = \begin{cases}
\frac{1}{2} & i \in \hat{C_r}, \, |\hat{C_r}|=2,\\
1 & i \in \hat{C_r}, \, |\hat{C_r}|=1.
\end{cases}\]
\end{defn}

\begin{thm}
\label{thm:NN_Greedy}
Let $\CC$ be a partial circular ordering ($|\CC|=m$) with a balanced TSP
weighting and $\delta$ a dissimilarity map.
A circular ordering $\CC'$ of size $|\CC'|=m-1$ that extends $\CC$ and
minimizes $l(\delta,\CC')$ is obtained by finding a pair $C_{r^*},C_{s^*}$ that minimize
\[ Q_{\delta}(C_r,C_s) = (m-2)\delta(C_r,C_s) - \sum_{t \neq r}\delta(C_r,C_t) - \sum_{t
  \neq s}\delta(C_s,C_t) \]
and then adding an edge between 
the pair of vertices corresponding to $i^* \in C_{r^*},j^* \in C_{s^*}$ in the circular ordering graph that minimize
\begin{eqnarray*}
 \hat{Q}_{\delta}(i,j) & = &(m-4+|\hat{C}_{r^*}|+|\hat{C}_{s^*}|)\delta(i,j) - \sum_{t \neq r^*,s^*}
\delta(i,C_t) - \sum_{t \neq r^*,s^*} \delta(j,C_t)\\ & &  - \sum_{k \in (C_{r^*} \cup C_{s^*})\setminus \{i\}} \delta(i,k) - \sum_{k \in (C_{r^*} \cup C_{s^*})\setminus \{j\}}
\delta(j,k). 
\end{eqnarray*}
\end{thm}

{\bf Proof}: Let $\CC=\{C_1,\ldots,C_m\}$ be a partial circular
ordering. A neighbor-net step consists of
adding an edge to $\CC$. This constitutes selecting two paths to join
(step 1), and then deciding which of the ends of the paths to join (step
2).
\begin{lemma}
The number of circular orderings consistent with $\CC$ is \[ |o(\CC)| = \frac{1}{2}(m-1)!
\prod_{r=1}^{m} |\hat{C}_r|.\]
\end{lemma}

\begin{figure}
\begin{picture}(500,500)
\put(20,-20){\includegraphics[scale=0.65]{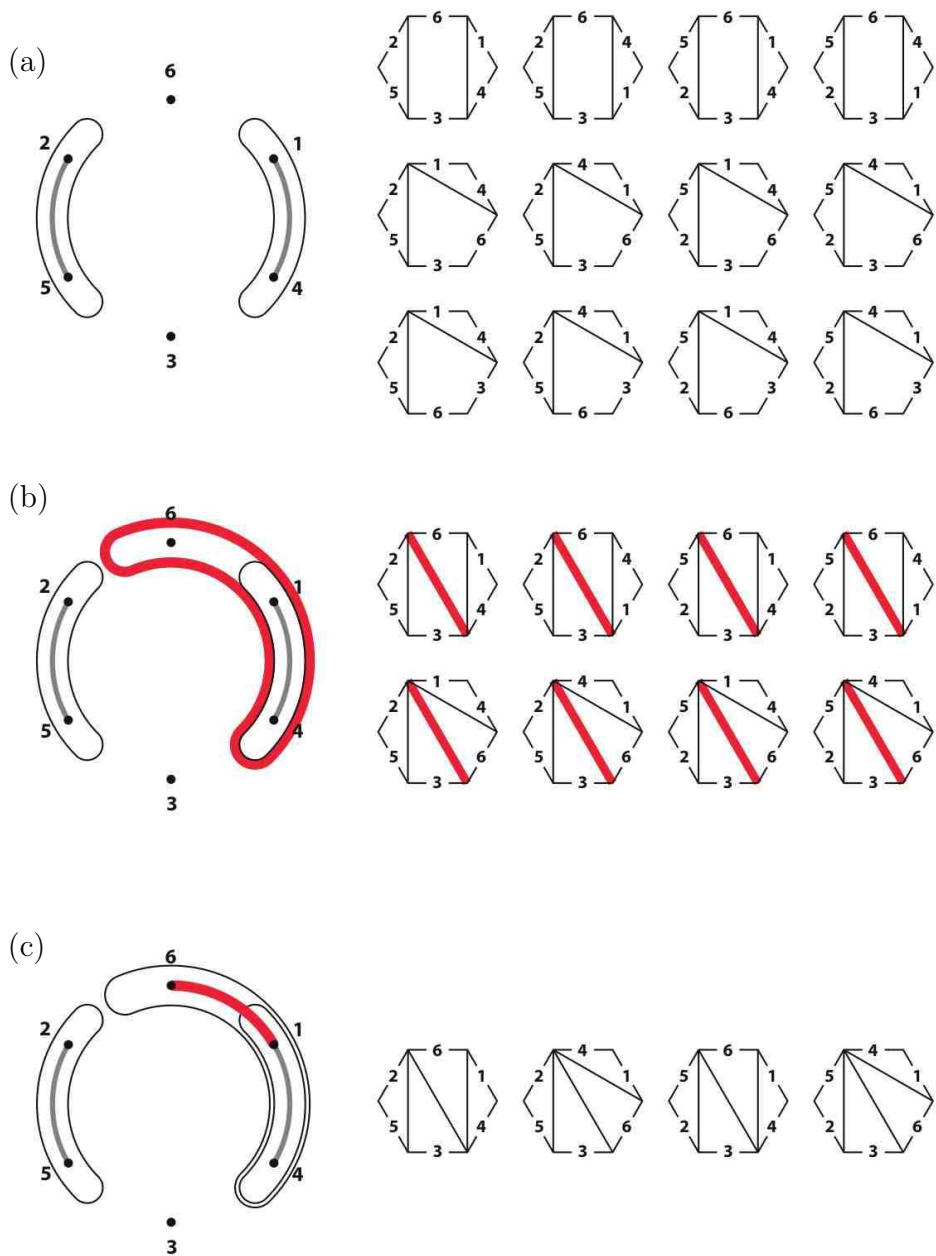}}
\put(40,465){(a)}
\put(40,300){(b)}
\put(40,130){(c)}
\end{picture}
\caption{A step in the neighbor-net algorithm run on the dissimilarity
  map $\delta$ from
  Figure 1(b). (a) A partial circular ordering $\CC, \, |\CC|=4$ and the 12 circular orderings
  consistent with it. Note that at this stage $l(\delta,\CC)=\frac{1708}{48}$ (b) Selection step part 1 showing $\CC_{r^*,s^*}$ where
$C_{r^*}=\{6\}$ and $C_{s^*}=\{1,4\}$. Now
  $l(\delta,\CC_{r^*,s^*})=\frac{1659}{48}$ and this is a neighbor-joining
agglomeration. (c) Selection step part 2
  results in a new partial circular ordering $\CC', \, |\CC'|=3|$ with 6
adjacent to 1 and $l(\delta,\CC')=\frac{1614}{48}$. This last step is what distinguishes
neighbor-net from neighbor-joining.}
\end{figure}

Let $\CC_{r,s}$ denote all of the partial circular orderings where there is
an edge between endpoints of $C_r$ and $C_s$ in the circular ordering
graph.
We say that a circular ordering is
consistent with $\CC_{r,s}$ if it is consistent with one of the partial
circular orderings in $\CC_{r,s}$. Similarly, we define $o(\CC_{r,s})$ to
constitute all circular orderings consistent with some partial circular ordering
 in $\CC_{r,s}$. In the following lemma we use the notation
 $\overline{ij}_{C}$ to denote that $i$ and $j$ are in the same
block in $C \in \CC$, and $i$ is adjacent to $j$.
\begin{lemma}
The number of circular orderings consistent with $\CC_{r,s}$ is
$2|o(\CC)|/(m-1)$ and
\[ \eta_{\CC_{r,s}}(i,j) =
 \begin{cases}
2|o(\CC)|/(m-1)& \mbox{ if } \overline{ij}_{C} \mbox{ for some } C,\\
2|o(\CC)|\mu(i)\mu(j)/(m-1) & \mbox{ if } i \in C_r, \, j \in C_s,\\
4|o(\CC)|\mu(i)\mu(j)/(m-1)(m-2)  & \mbox{ if } i \in C_t, \, j \in C_u, \, t \neq u, \, t,u \neq r,s,\\
2|o(\CC)|\mu(i)\mu(j)/(m-1)(m-2) & \mbox{ if } i \in C_r, \, j \in C_t, \, t
\neq s,\\
2|o(\CC)|\mu(x)\mu(y)/(m-1)(m-2) & \mbox{ if } i \in C_s, \,
j \in C_t, \, t \neq r,\\
0 & \mbox{ otherwise.}
\end{cases} \]
\end{lemma}
The proof of the lemma is elementary. We note that it also makes sense
for weightings that are not balanced TSP weightings, except that the
effect of the weightings $\mu$ is to alter the $\eta$ so that they count
the number of circular orderings consistent with split systems larger than
$\CC_{r,s}$.  For example, if $\mu$ is a tree weighting, then $\eta$ counts the
number of circular orderings consistent with the partially resolved tree $\mathcal{T}$.
For more on this see Definition \ref{def:length_splits} and
Theorem \ref{thm:variance}.

We can now conclude the proof of Theorem \ref{thm:NN_Greedy}:
\begin{eqnarray*}
\l(\delta,\CC_{r,s}) & = & \frac{1}{2}\sum_{C \in \CC}\sum_{\overline{ij}_C} \delta(i,j)
+ \frac{1}{2} \delta(C_{r},C_{s})
 + \frac{1}{(m-2)}\sum_{C_t \neq C_u, t,u \neq r,s}\delta(C_t,C_u)\\
& & + \frac{1}{2(m-2)} \sum_{C_t \neq C_{r},C_{s}}
 \delta(C_t,C_{s}) + \frac{1}{2(m-2)}\sum_{C_t \neq C_{r},C_{s}}  \delta(C_t,C_{r})\\
& = & \frac{1}{2}\sum_{C \in \CC}\sum_{\overline{ij}_C} \delta(i,j)+
 \frac{1}{(m-2)}\sum_{C_t \neq C_u}\delta(C_t,C_u) \\
& & + \frac{1}{2}\delta(C_{r},C_{s})
 - \frac{1}{2(m-2)} \sum_{C_t \neq C_{s}}
 \delta(C_t,C_{s}) - \frac{1}{2(m-2)}\sum_{C_t \neq C_{r^*}}  \delta(C_t,C_{s}).
\end{eqnarray*}

Thus, $l(D,\CC_{r,s}) = \frac{1}{2(m-2)}Q_{\delta}(C_r,C_s) + T$ where $T$ does not depend on $r$
or $s$. In other words, at each step neighbor-net is selecting a pair
$(r^*,s^*)$ to join that will minimize the balanced length. The actual
minimum balanced length is
attained for one of the $|\hat{C}_{r^*}||\hat{C}_{s^*}|$ possibilities for adding an edge
between $C_{r^*}$ and $C_{s^*}$ in $\CC$. Using the same argument as above, it
is easy to see that the minimum balanced length is attained when
$\hat{Q}_{\delta}(i,j)$ is subsequently minimized. \qed

\begin{remark}
  \label{rem:ZQ}
Let
\[ Z_{\delta}(C_r,C_s) =  -\frac{1}{m-1}\sum_{C_ \neq
  C_s}\delta(C_r,C_s)-\frac{1}{2}Q_{\delta}(C_r,C_s). \]
Then
\[ l(\delta,\CC) = \frac{1}{2}\sum_{C \in \CC}\sum_{i,j \in C}
\delta(i,j)+\frac{1}{(m-1)}T \]
implies that
\[ l(\delta,\CC) - l(\delta,\CC_{r,s}) = Z_{\delta}(C_r,C_s). \]
\end{remark}
The quantity $Z_{\delta}(C_r,C_s)$ features prominently in \cite{Contois2005,Gascuel1994,Mihaescu2006} and
is based on the ``neighborliness measurement'' of \cite{Gascuel1994}:
\[
Z_{\delta}(C_r,C_s) = \sum_{t,u \neq r,s} w(C_rC_s:C_tC_u), \mbox{ where }\]
\[ w(C_rC_s:C_tC_u)  =  \frac{1}{2}(\delta(C_r,C_t)
+ \delta(C_r,C_u) + \delta(C_s,C_t) + \delta(C_s,C_u) - 2 \delta(C_r,C_s) - 2\delta(C_t,C_u)).\]

It is interesting to note that the results in \cite{Mihaescu2006} are
motivated by this alternative formulation of the neighbor-joining
criterion. Remark \ref{rem:ZQ} provides further evidence that the
``$Z$-criterion'' is a natural formulation for the neighbor-joining
criterion, and at the same time explains the meaning of $Z_{\delta}(C_r,C_s)$ in
terms of the balanced length.

Returning to Remark \ref{rem:TSP}, we have the following interpretation of Theorem
\ref{thm:NN_Greedy}:
\begin{remark}
\label{rem:greedy}
Neighbor-net with a balanced TSP weighting is a greedy algorithm for the
traveling salesman problem.
\end{remark}
In fact, neighbor-net provides the optimal
solution for the TSP when $\delta$ satisfies the Kalmanson conditions (see Theorem
\ref{thm:Consistent} in Section 4). It is well known that the TSP can be
solved in polynomial time $O(n^2\mbox{log}n)$ for Kalmanson matrices \cite{Deineko1998}; neighbor-net
provides an alternative $O(n^3)$ polynomial algorithm. The $O(n^3)$
running time is based on the observation that the TSP and tree weighting
schemes can be implemented so that the selection steps are $O(k^2)$ where
$k$ is the number of blocks in the partial circular ordering at each
step. It should be
possible to obtain further improvements in speed by using the ideas developed for fast
neighbor-joining \cite{Elias2005}.

Theorem \ref{thm:NN_Greedy} is restricted to the balanced TSP
weighting.  We note, however, that there is no practical limitation to using different weightings
for the first and second selection steps.  We may consider a hybrid algorithm that applies a
tree weighting to the first selection step and a balanced TSP weighting to the second.  In that case,
Proposition \ref{prop:tree_construct} together with Theorem \ref{thm:NN_Greedy} show that
\begin{remark}
Neighbor-net with a hybrid weighting scheme is a greedy
algorithm for finding, simultaneously, the tree of minimum balanced length
 and the circular ordering of minimum length consistent with it.
\end{remark}

\section{The statistics}

We begin in this section by showing that neighbor-net is
a robust algorithm. By this we mean that if the input to neighbor-net is
a dissimilarity map  $\delta$ that is a perturbation of a circular decomposable
metric with respect to a circular ordering $\pi$, neighbor-net outputs
the circular ordering $\pi$. We note that in the case of a circular decomposable
metric where some of the splits
have zero weight, there will be more than one circular ordering consistent with $\delta$.  
In that case neighbor-net will output one of those circular orderings.  A corollary to this
is that if $\delta$ is a circular decomposable metric, and equation 
(\ref{eq:distance_estimates}) is used
to estimate the distances, then the output is exactly $\delta$, i.e.,
neighbor-net is a consistent estimator of the parameters of a circular
decomposable metric.  Implicit in the neighbor-net estimator are
assumptions about the variances of the measured distances. These can be
interpreted in terms of the weighting scheme used in neighbor-net, and we
return to this at the end of the section.

\begin{thm}
\label{thm:Consistent}
Suppose that $\delta:X \times X \rightarrow \mathbb{R}$ is a
dissimilarity map that satisfies the Kalmanson conditions for some
circular ordering $\pi$. Then neighbor-net applied to $\delta$ outputs a
circular ordering $\pi'$ such that $W_{\delta} \subseteq W_{\pi'}$.
\end{thm}
{\bf Proof}: It suffices to show that at any step of the algorithm, every circular ordering consistent with the partial circular ordering contains all the quartets in $W_{\delta}$.
Let $\CC=\{C_1,\ldots,C_m\}$ be a partial circular ordering consistent with
$\pi$ so that if $x_i \in C_r$ and $x_j \in C_s$ and $r<s$ then
$i<j$.
\begin{lemma}
\label{lem:Kalman}
For every $r < s < t < u$,
\begin{eqnarray*}
\delta(C_r,C_s) + \delta(C_t,C_u) & \leq & \delta(C_r,C_t) + \delta(C_s,C_u)\\
\delta(C_r,C_u) + \delta(C_s,C_t) & \leq & \delta(C_r,C_t) + \delta(C_s,C_u).
\end{eqnarray*}
\end{lemma}
{\bf Proof}: This follows directly from the Kalmanson conditions and the
requirement that $\sum_{i \in C_r}\mu(i) = 1$ for every $r$. \qed

Moreover, if for some $a \in C_r, b \in C_s, x \in C_t,y \in C_u$ with $\mu(a),\mu(b),\mu(x),\mu(y)>0$ we have $(ab;xy) \in W_{\delta}$, then $\delta(C_r,C_s)+\delta(C_t,C_u) < \delta(C_r,C_t)+\delta(C_s,C_u)$.

Next we introduce some notation to simplify the necessary calculations.
We set $\delta_{C_rC_s}(C_t) =
\delta(C_r,C_t)+\delta(C_s,C_t)-\delta(C_r,C_s)$. This is an analog of
the {\em Farris transform} \cite{Farris1977} for blocks in the partial circular ordering $\CC$.
Note that
\begin{equation}Q_{\delta}(C_r,C_s) = -2\delta(C_r,C_s) - \sum_{C_t}
  \delta_{C_rC_s}(C_t).\end{equation}
In order to simplify the presentation, we replace every $C_i$ with
$i$ in the formulas below. This is mathematically justified by Lemma
\ref{lem:Kalman} since blocks in a partial circular ordering
behave exactly like elements of the underlying set $X$ with respect
to the Kalmanson conditions. For example, by $Q_{\delta}(i,i+1)$ in the lemma below, we mean
$Q_{\delta}(C_i,C_{i+1})$ and a proof that $Q_{\delta}(C_i,C_{i+2})>Q_{\delta}(C_i,C_{i+1})$ is
equivalent to the proof that $Q_{\delta}(i,i+2)>Q_{\delta}(i,i+1)$ by Lemma \ref{lem:Kalman}.
\begin{lemma}
\label{lem:one_middle}
\[ Q_{\delta}(i,{i+2}) - Q_{\delta}(i,{i+1}) \geq 0. \]
\end{lemma}
{\bf Proof}: Let $j=i+2,\, k=i+1$.
\[ Q_{\delta}(i,j)-Q_{\delta}(i,k) = \sum_{x \neq i,j,k} \delta(k,x)+\delta(i,j) -
\delta(i,k) - \delta(j,x) \]
and $\delta(k,x)+\delta(i,j) -
\delta(i,k) - \delta(j,x) \geq 0$ for each $x$ by Lemma \ref{lem:Kalman}.
\qed
\begin{lemma}[The Anarchy Lemma]
\label{lem:two_middle}
\[ Q_{\delta}(i,{i+3})-Q_{\delta}({i+1},{i+2}) \geq 0. \]
\end{lemma}
{\bf Proof}: Let $j=i+3, \, k=i+1, \, l=i+2$. Applying Lemma
\ref{lem:Kalman} twice:
\begin{eqnarray*}
 Q_{\delta}(i,j)-Q_{\delta}(k,l) &  = & \sum_{x \neq i,j,k,l} (\delta(i,j) +
\delta(k,x)+\delta(l,x)) - \delta(i,x)-\delta(j,x)-\delta(k,l) \\
&\geq &  \sum_{x \neq i,j,k,l} (\delta(j,k) + \delta(i,x) +\delta(l,x) -\delta(i,x)-\delta(j,x)-\delta(k,l))\\
 &\geq &  \sum_{x \neq i,j,k,l}(\delta(j,k) +\delta(l,x) -\delta(j,x)-\delta(k,l)) \geq 0.\\
\end{eqnarray*} \qed

\begin{figure}
\includegraphics[scale=0.5]{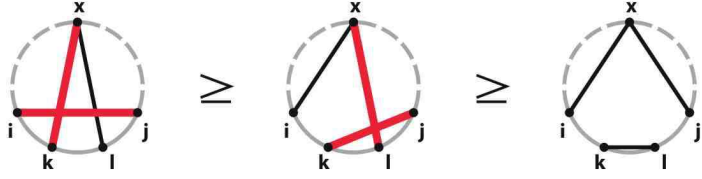}
\caption{An illustration of the proof of Lemma \ref{lem:two_middle}.}
\end{figure}

\begin{lemma}
\label{lem:3times}
Let $i<x<y<z<j<t$. Then
\begin{eqnarray*}
& &\delta_{xy}(z)+\delta_{xz}(y)+\delta_{yz}(x)+\delta_{xy}(t)+\delta_{xz}(t)+\delta_{yz}(t)\\
& \geq & 3\delta_{ij}(t)+\delta_{ij}(x)+\delta_{ij}(y)+\delta_{ij}(z).
\end{eqnarray*}
\end{lemma}
{\bf Proof}: Note that each of the following inequalities follows
directly from Lemma \ref{lem:Kalman}:
\begin{eqnarray*}
2\delta(x,t)+2\delta(i,j) & \geq &
\delta(x,i)+\delta(x,j)+\delta(i,t)+\delta(j,t)\\
2\delta(y,t)+2\delta(i,j) & \geq &
\delta(y,i)+\delta(y,j)+\delta(i,t)+\delta(j,t)\\
2\delta(z,t)+2\delta(i,j) & \geq &
\delta(z,i)+\delta(z,j)+\delta(i,t)+\delta(j,t).\\
\end{eqnarray*}
Summing both sides  we
obtain the required inequality. \qed
\begin{prop}
\label{prop:Consistent}
Suppose that $i < j - 3$. Then there exists $k$ such that
\begin{equation}
\label{eq:noncon}
 Q_{\delta}(i,j) - Q_{\delta}({k},{k+1})  \geq 0.
\end{equation}
\end{prop}

{\bf Proof}: Recall that $|\CC| = m$.  Suppose without loss of generality that $i = 0$ and $j \leq m/2$.
We will find $i < k \leq j-2$ satisfying (\ref{eq:noncon}), where the proof is
non-constructive and mimics the arguments in
Theorem 25 of \cite{Mihaescu2006}. In particular, we show that
\begin{equation}
\label{eq:main_nonconstr}
 (j - 3) \sum_{0 < x,y < j} \left( Q_{\delta}(i,j) - Q_{\delta}(x,y) \right) \geq 0,
\end{equation}
so that there exists $i<x,y<j$ with $Q_{\delta}(i,j)-Q_{\delta}(x,y) \geq 0$.

We first note that $$Q_{\delta}(i,j) - Q_{\delta}(x,y) = \sum_{t \neq i,j,x,y} \delta_{xy}(t) - \delta_{ij}(t).$$
We then break this sum into three sections: those that fall between $0$ and $j$, a matching set of the same size $(j-3)$
that lie beyond $j$, and lastly, all remaining terms.  In this way, equation (\ref{eq:main_nonconstr}) equals
$$(j-3)\sum_{0 < x, y <j}\left(\sum_{\substack{z = 1 \\ z \neq x, y}}^{j-1}\delta_{xy}(z) - \delta_{ij}(z)
                         + \sum_{t = j+1}^{2j-3}\delta_{xy}(t) - \delta_{ij}(t)
                         + \sum_{s = 2j-3}^{m-1}\delta_{xy}(s) - \delta_{ij}(s)\right)$$

By Lemma \ref{lem:two_middle}, the last summation is greater than or equal to zero, and so:

\begin{eqnarray*}
&\geq & \sum_{0 < x, y < j}(j-3)\left(\sum_{\substack{z = 1 \\ z \neq x, y}}^{j-1}\delta_{xy}(z) - \delta_{ij}(z)
                         + \sum_{t = j+1}^{2j-3}\delta_{xy}(t) - \delta_{ij}(t)\right)\\
& = & \sum_{0 < x, y < j} \quad  \sum_{\substack{z = 1 \\ z \neq x, y}}^{j-1} \quad \sum_{t=j+1}^{2j-3} \delta_{xy}(z)
+\delta_{xy}(t) - \delta_{ij}(z) - \delta_{ij}(t)\\
& = & \sum_{t=j+1}^{2j-3} \quad  \sum_{\substack{0 < x, y, z < j, \\ x \neq y \neq z}}
\delta_{xy}(z)+\delta_{xz}(y)+\delta_{yz}(x)
+\delta_{xy}(t)+\delta_{xz}(t)+\delta_{yz}(t)\\
& &
\qquad \qquad \qquad \qquad \qquad -3\delta_{ij}(t)-\delta_{ij}(x)-\delta_{ij}(y)-\delta_{ij}(z)
\quad \geq \quad 0.
\end{eqnarray*}
The final inequality follows from Lemma \ref{lem:3times}.  The claim
(\ref{eq:noncon}) now follows by noting that repeated application of the
argument leads to one of three cases: either we find a pair of neighbors
$k,k+1$ such that $Q_{\delta}(k,k+1) \leq Q_{\delta}(i,j)$, or else we
find a pair that are separated by one node (in which case we apply Lemma
\ref{lem:one_middle}) or a pair that are separated by two nodes (in
which case we apply Lemma \ref{lem:two_middle}). \qed

Returning to the proof of the theorem, it is clear that if 
we have a strict Kalmanson inequality on any quartet that separates $i$ and $j$, then 
the inequalities in Lemmas \ref{lem:one_middle}, \ref{lem:two_middle} and Proposition \ref{prop:Consistent} 
are strict inequalities.  Consequently we never join a pair of blocks that violate a quartet in $W_{\delta}$.
If the blocks are of size 1 we are done. Otherwise, it only remains to show that
two neighboring elements $x_r \in C_i$ and $x_{r+1} \in C_{i+1}$  will
be selected to be joined in the minimization of $\hat{Q}$. This follows
directly from the same arguments used in Lemmas \ref{lem:one_middle} and \ref{lem:two_middle}. \qed

The consistency of neighbor-net now follows easily by observing that for
a circular decomposable metric, the distances will be correctly
inferred using (\ref{eq:distance_estimates}).
\begin{cor}[\cite{Bryant2007}]
Neighbor-net is statistically consistent.
\end{cor}
Moreover, Theorem
\ref{thm:Consistent} can be used to obtain a neighbor-net analog of Atteson's theorem \cite{Atteson1999}
on the optimal radius of neighbor-joining:
\begin{cor}[Optimal radius]
\label{cor:optimal_radius}
Let $\mathcal{S}$ be a circular split system with respect to a circular
ordering $\pi=\{x_1,\ldots,x_n\}$, $\lambda_S > 0$ for
every $S \in \mathcal{S}$, and
$\delta_{\mathcal{S}}=\sum_{S \in \mathcal{S}} \lambda_S \delta_S$ a circular
decomposable metric. If $\epsilon = min_{S \in \mathcal{S}} \lambda_S$
and $\delta$ is any dissimilarity map with $||\delta-\delta_{\mathcal{S}}||_{\infty}
< \frac{\epsilon}{2}$ then neighbor-net will output a circular ordering whose split system contains $\mathcal{S}$.
\end{cor}
{\bf Proof}: It suffices to show that if $||\delta-\delta_{\mathcal{S}}||_{\infty}
\leq \frac{\epsilon}{2}$ then $\delta$ satisfies the Kalmanson
conditions with respect to $\pi$.
Let $i < j < k < l$.
\[
 \delta_{\mathcal{S}}(x_i,x_k)+\delta_{\mathcal{S}}(x_j,x_l) -
\delta_{\mathcal{S}}(x_i,x_j)-\delta_{\mathcal{S}}(x_k,x_l) = \sum_{S=\{A,B\},i,j\in A, k,l \in B}2\lambda_S.
\]
Therefore,
\[
 \delta(x_i,x_k)+\delta(x_j,x_l) -
\delta(x_i,x_j)-\delta(x_k,x_l) \geq \left( \sum_{S=\{A,B\},i,j\in A, k,l \in
  B}2\lambda_S \right) - 2 \epsilon > 0.
\]
A similar argument shows that  $\delta(x_i,x_k)+\delta(x_j,x_l) -
\delta(x_j,x_k)-\delta(x_k,x_i) \geq 0$.
\qed

Note that in Corollary \ref{cor:optimal_radius} the dissimilarity map $\delta$ satisfying $||\delta-\delta_{\mathcal{S}}||_{\infty} \leq \frac{\epsilon}{2}$ may not be a metric. Kalmanson matrices (as opposed to metrics) are characterized in \cite{Demidenko1997}.

We have already hinted at connections between neighbor-net and the
traveling salesman problem in Section 3. Our next theorem demonstrates
the consistency of the TSP estimate of the circular ordering and is
analogous to Theorem 2 of \cite{Desper2004}.
\begin{thm}
Let $\delta$ be a generic circular decomposable metric with respect to a circular
ordering $\pi=\{x_1,\ldots,x_n\}$. Then $l(\delta,\sigma)>l(\delta,\pi)$ for any circular
permutation $\sigma=\{y_1,\ldots,y_n\}$ different from $\pi$.
\end{thm}
{\bf Proof}: Since $\delta$ is a circular decomposable metric it must
satisfy the Kalmanson conditions. Therefore there must exist $i<k, \, |k-i|>1$
such that $\delta(y_i,y_{i+1})+\delta(y_k,y_{k+1}) >
\delta(y_i,y_k)+\delta(y_{i+1},y_{k+1})$. Consider the circular ordering

\[\sigma'= \{y_1,\ldots,y_i,y_{k},y_{k-1},\ldots,y_{i+1},y_{k+1},y_{k+2},\ldots,y_n\}. \]
Then $l(\delta,\sigma')<l(\delta,\sigma)$ and therefore
$\mbox{argmin}_{\tau} l(\delta,\tau) = \pi$. \qed

This result explains why it makes sense to use TSP solutions directly 
for finding circular orderings \cite{Korostensky2000}.

We now turn to the statistical meaning
of the weighting $\mu$ in the neighbor-net algorithm, and discuss how it should
be chosen in practice. We first consider the case of tree weightings.
In this case neighbor-net outputs a circular ordering consistent with
the neighbor-joining tree (Proposition \ref{prop:tree_construct}). The
theory of \cite{Desper2004} together with our results provides a direct
interpretation of the
agglomeration parameters that can be summarized as follows:
\begin{defn}[Length of a split system]
\label{def:length_splits}
Let $\mathcal{S}$ be a split system that is circular with respect
to some circular ordering and let $\eta_{\mathcal{S}}(i,j)$ be the number
of circular orderings
consistent with $\mathcal{S}$ where
where $x$ is adjacent to $y$. The length of a dissimilarity map $\delta$
with respect to $\mathcal{S}$ is
\[ l(\delta,\mathcal{S}) = \sum_{i,j}
\eta_{\mathcal{S}}(i,j)\delta(i,j).\]
\end{defn}
\begin{thm}
\label{thm:variance}
Let $\delta$ be a dissimilarity map, $\mathcal{S}$ a split system
that is circular with respect to some circular ordering, and
$\eta_{\mathcal{S}}(i,j)$ defined as above.
Let $\delta^*=\sum_{S \in \mathcal{S}}\lambda_S \delta_S$ ($\lambda_S
\geq 0$) be the circular decomposable metric obtained from the weighted least squares
estimates of the splits under the assumption that
the variance of $\delta(i,j)$ is
$\kappa \eta_{\mathcal{S}}(i,j)^{-1}$ (with the same constant $\kappa$
for all $i,j$). Then
\[ l(\delta,\mathcal{S}) = \sum_{S \in \mathcal{S}}\lambda_S. \]
\end{thm}
The choices of agglomeration
parameters for a tree weighting determine $\eta_{\mathcal{S}}(i,j)$ at
each step and are therefore implicit
variance assumptions
on the distances for the weighted least squares tree that is being
greedily approximated by the algorithm. The
balanced tree weighting scheme for neighbor-net corresponds to
balanced neighbor-joining agglomeration \cite{Desper2004}.
It should be interesting to
explore {\tt BIONJ}  \cite{Gascuel1997} analogs for neighbor-net,
which is easy to do since it only involves adapting the tree weightings.
In the case of a balanced TSP weighting, Theorem \ref{thm:variance}
explains that the neighbor-net algorithm ignores nodes once they have
two neighbors after agglomeration.

We conclude by remarking that some progress has been made in the
development of statistical models for split networks, suggesting the
possibility for maximum likelihood approaches to finding circular split
systems \cite{Bryant2005c,Sturmfels2007}.

\section{Applications}

Our goal in this section is to show how the theorems proved in the previous sections provide  insight into
how to use neighbor-net in practice, and in how to infer split networks.
We begin with an observation regarding the distance reduction formula
used in the current implementations of neighbor-net.

The agglomeration scheme proposed in \cite{Bryant2004} is as follows:
Suppose that a circular ordering contains two blocks $C_r,C_s$ that are
being agglomerated, where $C_r$ is a union of two smaller blocks
$C_r = C_t \cup C_u$ so that the
agglomerated block is $C_r \cup C_s = C_t \cup C_u \cup C_s$ in
that order.
\begin{equation}
\label{eq:oldagglom}
\mu'(i) = \begin{cases}
\frac{1}{4}\mu(i) & i \in C_t \cup C_s, \\
\frac{1}{2}\mu(i) & i \in C_u.
\end{cases}\end{equation}
There is an analogous formula for the case when two blocks, each
composed of two blocks are being joined (the above formula is applied
twice).

This weighting is neither a TSP weighting nor a tree
weighting. Furthermore, in the case of agglomeration of a pair of blocks
each composed of two blocks, the resulting weighting $\mu$ depends on the
order in which the agglomeration is performed. Thus,
the tree output by neighbor-net using (\ref{eq:oldagglom}) is not
necessarily the neighbor-joining tree, whereas the use of a
tree-weighting scheme guarantees this (Proposition
\ref{prop:tree_construct}).

The advantage of producing a circular ordering consistent with the
neighbor-joining tree, is that it allows for a direct analysis of the
conflicting signals with a tree of interest.
To demonstrate this, we analyzed a published dataset of language
structure characters from Oceanic Austronesian and Papuan languages
\cite{Dunn2005}. The neighbor-net algorithm was previously used to
infer phylogenetic  relationships among the languages (Figure S2 from
the supplementary materials of \cite{Dunn2005}).
We compared Figure S2 obtained using the default parameters for
neighbor-net (\ref{eq:oldagglom}) with the balanced tree weighting scheme that produces a
neighbor-joining tree.  In both cases,
the split weights were computed using the constrained least squares
estimation procedure in \cite{Bryant2004}. The split networks were
visualized using the program {\tt SplitsTree4} \cite{Huson2005}. Figure 4(a)
shows the network for the balanced tree weighting scheme, together with
the neighbor-joining tree corresponding to the split system output by
the algorithm. The circular ordering obtained by using the default
neighbor-net settings
is not consistent with this neighbor-joining
tree. The ability to view the neighbor-joining tree in conjunction with
the neighbor-net split network is a direct result of Proposition \ref{prop:tree_construct}.
The representation of the tree together with the
network, as shown in Figure 4(left), is useful for directly using
 neighbor-net to evaluate the extent of phylogenetic discordancy with the
 neighbor-joining tree. For example, we see clearly that the split
 between the Papuan and Austronesian (Oceanic) languages is in fact a
 split in the neighbor-joining tree. Note that all the edges in the
 network and tree are drawn to scale.

\begin{figure}
\includegraphics[scale=0.315]{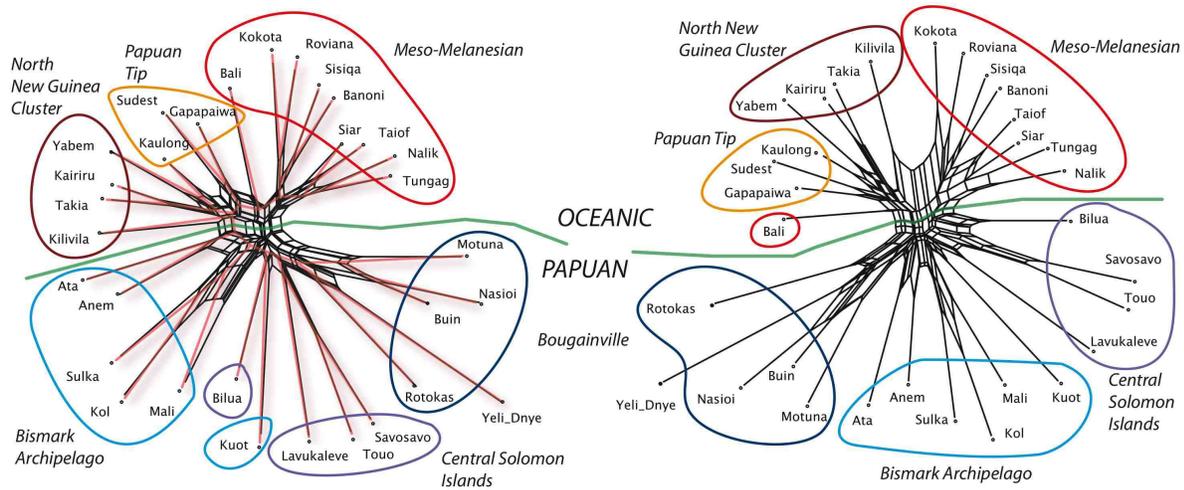}
\caption{Left: Neighbor-net and the neighbor joining tree for groups of
  Papuan and Austronesian languages. Right: The
  split-network inferred for the optimal circular ordering that was
  obtained using {\tt Concorde}.}
\end{figure}

The interpretation of neighbor-net as a greedy algorithm for the TSP
suggests an analysis of the optimal TSP tour. We computed this tour for
the dataset from \cite{Dunn2005} using {\tt Concorde} \cite{Concorde}. The optimal tour, of
length 7.541 was found in 0.57 seconds. The length of this tour should be
contrasted with the length of the balanced tree weighting tour, 7.810,
which is very close to 7.794, the length of the tour obtained using the default
parameters. The constrained least squares optimization procedure
of \cite{Bryant2004} was applied to the optimal circular ordering and
resulted in the split network shown in Figure 4(right).

The comparison of the two split networks in Figure 4 is interesting. A
key observation in \cite{Dunn2005} was that the Papuan languages cluster
into groups consistent with the geographical locations of the islands.
On the other hand, it was remarked that Bougainville, which is
geographically in between the Bismarck Archipelago and the Central
Solomon Islands, did not cluster in between the languages from those two
locations. Figure 4(right) shows that the TSP ordering produces a better
overall clustering, albeit still with the Bougainville languages not
sandwiched in the geographically correct location. Nevertheless, a key
new insight that emerges from the network is that Bali, which appears to
be incorrectly grouped, is in fact correctly grouped if one assumes that
the Papuan and Oceanic groups are really two distinct separate groups
(it is then just a neighbor to Nalik).

Our main conclusion is that the
choice of weightings in the neighbor-net algorithm is important in
determining the results, and that care has to be taken in choosing the
weights appropriately. Furthermore, tree weighting algorithms will be
useful in cases where it is desirable to use neighbor-net as a
diagnostic tool for exploring neighbor-joining trees, and TSP algorithms
may be useful for direct application in obtaining circular
orderings. In fact, the use of TSP solvers in
similar contexts is not new, appearing in \cite{Korostensky2000} in
the context of tree construction and in \cite{Johnson2006}, where the  {\tt Concorde} program  is used
to find a circular ordering from a distance matrix for proteins based on
protein-protein interactions.
It also seems important to develop a variant of neighbor-net
that outputs the optimal circular ordering consistent with an arbitrary given tree.

We conclude by noting that neighbor-net can also be used practically as a greedy
algorithm for the TSP. Unlike the naive greedy algorithm for which
many negative results have been published (see, e.g., \cite{bangjensen2004}),
neighbor-net exhibits good properties. For example, the output does not
depend on the order of the input, and the algorithm is optimal for Kalmanson
matrices.
We experimented with the problem st70.tsp from {\tt TSPLIB} \cite{Reinelt1991}. The
balanced TSP weighting gave a tour of length 759.801, that is only 12\%
longer than the optimal tour of length 678.598. As expected, the
balanced tree
weighting scheme yielded a longer tour of length 812.613. It will be
interesting to explore the improvements possible with the incorporation
of search heuristics such as nearest
neighbor interchange moves. These have been used to
significantly improve neighbor-joining in the {\tt FastME} program
 \cite{Desper2002}.

\section*{Acknowledgments}
Figure 1(e) is based on Figure 13(b) from \cite{Devadoss2004} and we thank
Stayan Devadoss for permission to use his artwork. We thank David
Bryant, Vincent Moulton and Andreas Spillner for kindly sharing with us a preprint of \cite{Bryant2007}, Lu Luo for useful comments on a preliminary version of this manuscript, and Radu
Mihaescu for suggestions in the proof of Theorem \ref{thm:Consistent}. 
Lior Pachter was supported in part by an NSF CAREER award (CCF-0347992)
and thanks Jotun Hein and Philip Maini for hosting him while on
sabbatical when this work was performed.
 Dan Levy was supported by a grant from the
Biotechnology and Biological Sciences Research Council of the UK
(BB/D005418/1).

\newpage

\bibliographystyle{amsplain}
%\bibliography{zproof}

\begin{thebibliography}{10}

\bibitem{Concorde}
D~Applegate, R~Bixby, V~Chvatal, and W~Cook, \emph{{The Concorde TSP solver}},
  {\tt http://www.tsp.gatech.edu/concorde.html/}.

\bibitem{Atteson1999}
K~Atteson, \emph{The performance of neighbor-joining methods of phylogenetic
  reconstruction}, Algorithmica \textbf{25} (1999), 251--278.

\bibitem{Bandelt1992}
HJ~Bandelt and A~Dress, \emph{A canonical decomposition theory for metrics on a
  finite set}, Advances in Mathematics \textbf{92} (1992), 47--105.

\bibitem{bangjensen2004}
J~Bang-Jensen, G~Gutin, and A~Yeo, \emph{When the greedy algorithm fails},
  Discrete Optimization \textbf{1} (2004), 121--127.

\bibitem{Barad2003}
G~Barad, \emph{Genome rearrangements and algebraic geometry}, Knots in
  Washington XV (K~Kobayashi, K~Przytycki, Y~Rong, S~Suzuki, K~Taniyama,
  T~Tsukamoto, and A~Yasuhara, eds.), 2003.

\bibitem{Billera2001}
LJ~Billera, SP~Holmes, and K~Vogtmann, \emph{Geometry of the space of
  phylogenetic trees}, Advances in Applied Mathematics \textbf{27} (2001),
  733--767.

\bibitem{Bryant2005c}
D~Bryant, \emph{Extending tree models to split networks}, Algebraic Statistics
  for Computational Biology (L~Pachter and B~Sturmfels, eds.), Cambridge
  University Press, 2005, pp.~322--334.

\bibitem{Bryant2005b}
D~Bryant, F~Filimon, and R~Gray, \emph{Untangling our past: languages, trees,
  splits and networks}, The evolution of cultural diversity: phylogenetic
  approaches (R~Mace, C~Holden, and S~Shennan, eds.), UCL Press, 2005,
  pp.~69--85.

\bibitem{Bryant2007}
D~Bryant, V~Moulton, and A~Spillner \emph{Consistency of the {N}eighbor{N}et algorithm for
  constructing phylogenetic networks}, Algorithms for Molecular Biology \textbf{2} (2007).

\bibitem{Bryant2004}
\bysame, \emph{{NeighborNet}: An agglomerative method for the construction
  of planar phylogenetic networks}, Molecular Biology And Evolution \textbf{21}
  (2004), 255--265.

\bibitem{Buneman1971}
P~Buneman, \emph{The recovery of trees from measures of dissimilarity},
  Mathematics in the Archaeological and Historical Sciences (FR~Hodson,
  DG~Kendall, and P~Tautu, eds.), Edinburgh University Press, 1971,
  pp.~387--395.

\bibitem{Carr2006}
M~Carr and S~Devadoss, \emph{Coxeter complexes and graph associahedra},
  Topology and its applications \textbf{153} (2006), 2155--2168.

\bibitem{Chepoi1998}
V~Chepoi and B~Fichet, \emph{A note on circular decomposable metrics},
  Geometrica Dedicata \textbf{69} (1998), 237--240.

\bibitem{Christopher1997}
G~Christopher, \emph{Structure and applications of totally decomposable
  metrics}, Ph.D. thesis, Carnegie Mellon University, 1997.

\bibitem{Christopher1996}
G~Christopher, M~Farach, and M~Trick, \emph{The structure of circular
  decomposable metrics}, Lecture Notes in Computer Science, vol. 1136,
  Springer, New York, 1996, pp.~406--418.

\bibitem{Contois2005}
M~Contois and D~Levy, \emph{Small trees and generalized neighbor-joining},
  Algebraic Statistics for Computational Biology (L~Pachter and B~Sturmfels,
  eds.), Cambridge University Press, 2005, pp.~333--344.

\bibitem{Deineko1998}
VG~Deineko, R~Rudolf, and GJ~Woeginger, \emph{Sometimes traveling is easy: the
  master tour problem}, SIAM Journal of Discrete Mathematics \textbf{11}
  (1998), 81--93.

\bibitem{Demidenko1997}
VM~Demidenko and R~Rudolf, \emph{A note on Kalmanson matrices}, Optimization \textbf{40} (1997), 285--294.

\bibitem{Desper2002}
R~Desper and O~Gascuel, \emph{Fast and accurate phylogeny reconstruction
  algorithms based on the minimum-evolution principle}, Journal of
  Computational Biology \textbf{19} (2002), no.~5, 687--705.

\bibitem{Desper2004}
\bysame, \emph{Theoretical foundation of the balanced minimum evolution method
  of phylogenetic inference and its relationship to weighted least-squares tree
  fitting}, Molecular Biology and Evolution \textbf{21} (2004), 587--598.

\bibitem{Desper2005}
\bysame, \emph{The minimum evolution distance-based approach to phylogenetic
  inference}, Mathematics of Evolution and Phylogeny (O~Gascuel, ed.), Oxford
  University Press, 2005.

\bibitem{Devadoss1999}
S~Devadoss, \emph{Tessellations of moduli spaces and the mosaic operad},
  Contemporary mathematics \textbf{239} (1999), 91--114.

\bibitem{Devadoss2004}
\bysame, \emph{Combinatorial equivalence of real moduli spaces}, Notices of the
  American Mathematical Society \textbf{51} (2004), 620--628.

\bibitem{Dress2004}
A~Dress and DH~Huson, \emph{Constructing splits graphs}, IEEE/ACM Transactions
  in Computational Biology and Bioinformatics \textbf{1} (2004), 109--115.

\bibitem{Dress1996}
A~Dress, V~Moulton, and W~Terhalle, \emph{T-theory: an overview}, European
  Journal Combinatorics \textbf{17} (1996), 161--175.

\bibitem{Dunn2005}
M~Dunn, A~Terrill, G~Reesnik, RA~Foley, and SC~Levinson, \emph{Structural
  phylogenetics and reconstruction of ancient language history}, Science
  (2005), 2072--2075.

\bibitem{Elias2005}
I~Elias and J~Lagergren, \emph{Fast neighbor joining}, Proceedings of the
  International Colloquium on Automata, Languages and Programming (ICALP '05),
  2005.

\bibitem{Farris1977}
JS~Farris, \emph{On the phenetic approach to vertebrate classification}, Major
  patterns in vertebrate evolution, Plenum, New York, 1977, pp.~823--950.

\bibitem{Gascuel1994}
O~Gascuel, \emph{{A note on Sattath and Tversky's, Saitou and Nei's, and
  Studier and Keppler's Algorithms for Inferring Phylogenies from Evolutionary
  Distances}}, Molecular Biology and Evolution \textbf{11} (1994), 961--963.

\bibitem{Gascuel1997}
\bysame, \emph{{BIONJ: an improved version of the NJ algorithm based on a
  simple model of sequence data}}, Molecular Biology and Evolution \textbf{14}
  (1997), 685--695.

\bibitem{Gascuel2006}
O~Gascuel and M~Steel, \emph{{Neighbor-joining revealed}}, Molecular Biology
  and Evolution \textbf{23} (2006), 1997--2000.

\bibitem{Huson2005b}
D~Huson and D~Bryant, \emph{Application of phylogenetic networks in
  evolutionary studies}, Molecular Biology and Evolution \textbf{23} (2005),
  254--267.

\bibitem{Huson2005}
\bysame, \emph{{Estimating phylogenetic trees and networks using SplitsTree4}},
  in preparation, 2005.

\bibitem{Johnson2006}
O~Johnson and J~Liu, \emph{A traveling salesman approach for predicting protein
  functions}, Source Code for Biology and Medicine \textbf{1} (2006).

\bibitem{Kalmanson1975}
K~Kalmanson, \emph{Edgeconvex circuits and the traveling salesman problem},
  Canadian Journal of Mathematics \textbf{27} (1974), 1000--1010.

\bibitem{Korostensky2000}
C~Korostensky and G~Gonnet, \emph{Using traveling salesman problem algorithms
  for evolutionary tree construction}, Bioinformatics \textbf{16} (2000),
  619--627.

\bibitem{Lutz2006}
M~Lutz, \emph{Bantu classification, {B}antu trees and phylogenetic methods},
  Phylogenetic Methods and the Prehistory of Languages (P~Foster and C~Renfrew,
  eds.), Cambridge: McDonald Institute for Archaeological Research, 2006,
  pp.~43--55.

\bibitem{Mihaescu2006}
R~Mihaescu, D~Levy, and L~Pachter, \emph{Why neighbor joining works}, arXiv
  cs.DS/0602041, 2006.

\bibitem{Moret2004}
BE~Moret, L~Nakhleh, T~Warnow, CR~Linder, A~Tholse, A~Padolina, J~Sun, and
  R~Timme, \emph{Phylogenetic networks: modeling reconstructibility and
  accuracy}, IEEE/ACM Transactions on Computational Biology and Bioinformatics
  \textbf{1} (2004), 13--23.

\bibitem{Newberg1996}
LA~Newberg, \emph{The number of clone orderings}, Discrete Applied Mathematics
  \textbf{69} (1996), 233--245.

\bibitem{ASCB2005}
L~Pachter and B~Sturmfels (eds.), \emph{Algebraic statistics for computational
  biology}, Cambridge University Press, 2005.

\bibitem{Reinelt1991}
G~Reinelt, \emph{{TSPLIB - A traveling salesman problem library}}, ORSA Journal
  on Computing \textbf{3} (1991), 376--384.

\bibitem{Saitou1987}
N~Saitou and M~Nei, \emph{The neighbor joining method: a new method for
  reconstructing phylogenetic trees}, Molecular Biology and Evolution
  \textbf{4} (1987), 406--425.

\bibitem{Semple2003}
C~Semple and M~Steel, \emph{Phylogenetics}, {\rm Oxford Lecture Series in
  Mathematics and its Applications}, vol.~24, Oxford University Press, Oxford,
  2003.

\bibitem{Semple2004}
\bysame, \emph{Cyclic permutations and evolutionary trees}, Advances in Applied
  Mathematics \textbf{32} (2004), 669--680.

\bibitem{Stanley1999}
RP~Stanley, \emph{Enumerative combinatorics. {V}ol. 2}, Cambridge Studies in
  Advanced Mathematics, vol.~62, Cambridge University Press, Cambridge, 1999.

\bibitem{Sturmfels2007}
B~Sturmfels and S~Sullivant, \emph{Toric geometry of cuts and splits}, arXiv
  math.AC/0606683, 2007.

\end{thebibliography}
\providecommand{\bysame}{\leavevmode\hbox to3em{\hrulefill}\thinspace}
\providecommand{\MR}{\relax\ifhmode\unskip\space\fi MR }
\providecommand{\MRhref}[2]{%
  \href{http://www.ams.org/mathscinet-getitem?mr=#1}{#2}
}
\providecommand{\href}[2]{#2}

\end{document}